\documentclass{amsart}
\usepackage[utf8]{inputenc}

\usepackage{setspace}
\usepackage[a4paper]{geometry}
\usepackage{graphicx}
\DeclareGraphicsExtensions{.pdf,.png,.jpg}

\usepackage{caption}

\usepackage{enumerate}
\usepackage{float}
\usepackage{todonotes}

\usepackage{hyperref}
\usepackage{cleveref}

\usepackage{dsfont}
\usepackage[mathscr]{euscript}
\usepackage{mathrsfs}
\usepackage{mathtools}

\usepackage{tikz-cd}

\newcommand{\calF}{\mathcal{F}}

\newcommand{\N}{\mathbb{N}}

\newcommand{\vp}{\varphi}

\renewcommand{\vp}{\varphi}

\newcommand{\blank}{\underline{\hspace{.3 cm}}}

\usepackage[symbol]{footmisc}

\input{xy}
\xyoption{all}
\usepackage{tikz}

\newtheorem{thm}{Theorem}[section]

\newtheorem{prop}[thm]{Proposition}

\newtheorem{question}[thm]{Question}

\theoremstyle{definition}
\newtheorem*{dfn}{Definition}

\theoremstyle{remark}
\newtheorem*{remark}{Remark}

\usepackage{amsthm,amsmath,amstext,amsfonts,amssymb}
\makeatother

\DeclareMathOperator{\Hom}{Hom}
\DeclareMathOperator{\Ext}{Ext}

\makeatletter
\newenvironment{myalign*}{\ifvmode\else\hfil\null\linebreak\fi
  \hspace*{-\leftmargin}\minipage\textwidth
  \setlength{\abovedisplayskip}{0pt}%
  \setlength{\abovedisplayshortskip}{\abovedisplayskip}%
  \start@align\@ne\st@rredtrue\m@ne}%
{\endalign\endminipage\linebreak}
\makeatother

\makeatletter
\@namedef{subjclassname@2020}{%
  \textup{2020} Mathematics Subject Classification}
\makeatother

\title{A remark on the category of graded $F$-modules}
\author{McKinley Gray}
\address{}
\email{}
\subjclass[2020]{}

\begin{document}

\maketitle

\begin{abstract}
Let $R=k[x,y]$ be a polynomial ring over a field $k$ of prime characteristic $p$ and let $E$ denote the injective hull of $k$ (which is isomorphic to $H^2_{(x,y)}(R)$). We prove that $E$ is not an injective object in the category of graded $F$-modules over $R$. This answers in the negative a question raised by Lyubeznik-Singh-Walther.
\end{abstract}

\section{Introduction}
Building upon previous research such as \cite{PeskineSzpiro, HartshorneSpeiser, HunekeSharp}, Lyubeznik introduces in \cite{LyuFModule} the theory of $F$-modules over a noetherian regular ring of prime characteristic $p$. We will review this theory in the next section. For the purpose of the introduction, we remark that local cohomology modules $H^j_I(R)$ supported in an ideal $I$ of a noetherian regular ring $R$ of prime characteristic are primary examples of ($F$-finite) $F$-module. Lyubeznik's theory of $F$-modules has become indispensable in the study of rings of prime characteristic $p$, see for instance \cite{EmertonKisin, WZhang, LSW}. When $R=k[x_1,\dots,x_n]$ is a polynomial ring over a field $k$ of prime characteristic with its natural grading, one may consider the notion of graded $F$-module ({\it cf.} \cite{YiZhang}). In the theory of graded $F$-modules, one particular local cohomology module plays a special role, namely $H^n_{(x_1,\dots,x_n)}(R)$, which we will denote by $E$ throughout this paper, as it is isomorphic to the (graded) injective hull of $k$ as an $R$-module. 

In \cite{LSW}, the authors prove that $E$ is an injective object in the category of {\it graded} $F$-finite $F$-modules, which is surprising since it is proved in \cite{Ma14} that $E$ is not an injective object in the category of $F$-finite $F$-modules. This prompts the question on the injectivity of $E$ in the category of graded $F$-modules. Indeed, the following is left open in \cite[2.11]{LSW}.

\begin{question}
\label{main question}
Let $R=k[x_1,\dots,x_n]$ be a standard graded polynomial ring and $E$ be as above. Is $E$ an injective object in the category of graded $F$-modules?
\end{question}

The main purpose of this paper is to answer Question \ref{main question} in the negative by proving the following.

\begin{thm}
\label{main theorem}
Let $R=k[x,y]$ be a standard graded polynomial ring over a perfect field $k$ of characteristic $p$ and $M=(R_x)^\vee$ (the graded Matlis dual of $R_x$). Then, 
\[\Ext^1_{\calF}(M,E)\neq 0\]
where $\calF$ denotes the category of graded $F$-modules, and $E=H^2_{(x,y)}(R)$. Consequently, $E$ is not an injective object in the category of graded $F$-modules.
\end{thm}

Theorem \ref{main theorem} is proved, in \S\ref{proof of main theorem}, by constructing a non-splitting short exact sequence in the category of graded $F$-modules:
\[0\to E\to L\to M\to 0.\]

In \S\ref{prelim}, we will review some basic facts of (graded) $F$-modules.

\subsection*{Acknowledgments}
The author would like to thank Wenliang Zhang, his thesis advisor, for patient guidance and support and many valuable discussions during the work on this result.  The author would also like to thank the anonymous referee whose constructive feedback improved the quality of this paper.
\section{Preliminaries}
\label{prelim}
\begin{dfn}
Let $R$ be a noetherian regular ring of prime characteristic $p$ and let $F:R\to R$ be the Frobenius endomorphism. We will denote the functor $F_*R\otimes_R -$ by $F(-)$. An $R$-module $M$ is an {\it $F_R$-module} (or {\it $F$-module} when $R$ is understood) if there is an $R$-module isomorphism $M\xrightarrow{\sim}F(M)$. An $F$-module morphism between two $F$-modules $M,N$ is a commutative diagram:

\[
\xymatrix{
M \ar[r]^{\vp} \ar[d] & N \ar[d]\\
F(M) \ar[r]^{F(\vp)} & F(N)
}
\]
An $F$-module $M$ is {\it $F$-finite} if there is a finitely generated $R$-module $\tilde{M}$ with an $R$-module homomorphism $\beta:\tilde{M}\to F(\tilde{M})$ such that $M$ is isomorphic to the direct limit $\varinjlim(\tilde{M}\xrightarrow{\beta} F(\tilde{M})\xrightarrow{F(\beta)} F^2(\tilde{M})\cdots)$.
\end{dfn}

\begin{dfn}
A {\it graded $F$-module} is a graded $R$-module $M$ with a degree-preserving $R$-module isomorphism $M\xrightarrow{\sim}F(M)$ with the grading on $F(M)$ given by $\deg(r\otimes m)= \deg(r)+ p \deg(m) $

Graded $F$-finite $F$-modules can be defined analogously. 
\end{dfn}

\begin{prop}\label{prop1} For a given graded $F$-module, $N$, with structure morphism $\theta_N$, applying ${}^{*}\Hom(\blank,E)$ to $\theta_N^{-1}$ produces a degree preserving isomorphism ${}^{*}\Hom(N,E)\rightarrow {}^{*}\Hom(F(N),E)$. 
\end{prop}
\begin{proof}
$\theta_N^{-1}$ is degree preserving, so pre-composing by it will not change the degree of a map to $E$.  $E$ is injective as an $R$-module and $\theta_N^{-1}$ is an isomorphism so the induced map is an isomorphism.  
\end{proof}

\begin{prop}\label{prop2}
There is a degree preserving isomorphism \[\sideset{^*}{_R}\Hom(F(R_x),E)\xrightarrow{\sim} F(\sideset{^*}{_R}\Hom(R_x,E))\]
\end{prop}

\begin{proof} To get this morphism we first take $f\in {}^*\Hom(F(R_x),E)$ of some given degree $d$.  Find $\alpha$ and $\beta$ such that $d=(\alpha-1)p-\beta$ with $-p+2\leq \beta\leq 1$ so that we can write $f$'s output of all $1\otimes \frac{1}{x^\ell}$ as follows

\begin{align*}
    f\left(1\otimes \frac{1}{x^{\alpha-1}}\right)&=0\\
f\left(1\otimes \frac{1}{x^\alpha}\right)= \frac{c_0}{x^{p}y^\beta}+ \frac{c_1}{x^{p-1}y^{\beta+1}}+&\ldots+\frac{c_{p-1}}{xy^{\beta+p-1} }  \\
f\left(1\otimes \frac{1}{x^{\alpha+1}}\right)= \frac{c_0}{x^{2p}y^\beta}+ \frac{c_1}{x^{2p-1}y^{\beta+1}}+&\ldots +\frac{c_{p-1}}{x^{p+1}y^{\beta+p-1} }+ \frac{c_p}{x^{p}y^{\beta+p}}+\ldots + \frac{c_{2p-1}}{xy^{\beta+2p-1}} \\
\vdots\\
f\left(1\otimes \frac{1}{x^{\alpha+n}}\right)&=\sum\limits_{i=0}^{(n+1)p-1} \frac{c_i}{x^{(n+1)p-i}y^{\beta+i}}
\intertext{Note that many of these terms may in fact be zero, namely the first $1-\beta$ terms will be zero because of the $y$-degree. Furthermore, many $c_i$ are potentially zero, thus $f(1\otimes \frac{1}{x^\alpha})$ may be zero and thus the first nonzero term would have a larger power of $y$.  The above formatting is as generic as possible for an $f$ with the given degree.  }
\intertext{Define  a collection of maps as follows}
f_0\left(1\otimes \frac{1}{x^{\alpha+n}}\right)&:=\frac{c_0}{x^{(n+1)p} y^\beta} \ \ \ \ \ \ \ \forall n\geq 0\\
f_1\left(1\otimes \frac{1}{x^{\alpha+n}}\right)&:=\frac{c_1}{x^{(n+1)p-1} y^{\beta+1}} \ \ \ \ \ \ \ \forall n\geq 0\\
\vdots\\
\intertext{for a general $i=qp+r$ with $0\leq r<p-1$ define }
f_i\left(1\otimes \frac{1}{x^{\alpha+n}}\right)&= \frac{c_i}{x^{(n+1)p -i}y^{\beta+i}}\ \ \ \ \ \forall n\geq q
\intertext{That is to say, each $f_i$ just outputs one specific term of $f$ (and continues to only output one term as the degree of the input descends), so that we can formally rewrite $f$ as }
f&=\sum\limits_{i\in \N} f_i\\
\intertext{The $f_i$'s have associated $m_i, n_i, k_i ,\ell_i $ with $0\leq m_i, n_i<p$ for which}
f_i\left(x^{m_i}y^{n_i}\otimes \frac{y^{k_i}}{x^{\ell_i}}\right)&=\frac{c_i}{xy}\\
\intertext{So $(\ell_i-k_i)p-m_i-n_i-2=d$.  Define $\hat{f}_i\in {}^*\Hom(R_x,E)$ as follows }
\hat{f}_i\left(\frac{1}{x^{\ell_i+n}}\right)&=\frac{(c_i)^{\frac{1}{p}}}{x^{n+1}y^{k_i+1}}\ \ \ \ \forall n\in \N\\
\end{align*}

So $\hat{f}_i$ have degree $\ell_i-k_i-2$.

Now we define $\sideset{^*}{_R}\Hom(F(R_x),E)\xrightarrow{\sim} F(\sideset{^*}{_R}\Hom(R_x,E))$

\[f\mapsto \sum\limits_{i\in \N} \left(x^{p-m_i-1}y^{p-n_i-1}\otimes \hat{f}_i\right)\]
perhaps appropriately written to see that this infinite sum can be made sense of in this case as

\[f\mapsto \sum\limits_{0\leq n<p} \left( x^{p-m_i-1}y^{p-n_i-1}\otimes \sum\limits_{i\equiv n\  \text{mod}\ p} \hat{f}_i \right)\]

Terms have degree $(\ell_i-k_i-2)p+2p-m_i-n_i-2= (\ell_i-k_i)p-m_i-n_i-2=d$ so the map is degree preserving.

To check that this is an isomorphism we can check that this map is the inverse to

\[r\otimes f\mapsto \left(\left(1\otimes \frac{1}{x^\ell}\right)\mapsto r f\left(\frac{1}{x^\ell}\right)^p\right)
\]
 
Take $x^a y^b\otimes g\in F_*R \otimes \Hom(R_x,E)$ with $0 \leq a,b< p$, $g$ of degree $\alpha-\beta-1$ with output written as follows (again with many $c_i$ potentially zero) 

\begin{align*}
    g\left(\frac{1}{x^{\alpha-1}}\right)&=0\\
    g\left(\frac{1}{x^\alpha}\right)&= \frac{c_0}{xy^\beta}\\
    &\vdots\\
    g\left(\frac{1}{x^{\alpha+n}}\right)&= \sum\limits_{i=0}^{n} \frac{c_i}{x^{n-i+1}y^{\beta+i}}\\ 
    \intertext{Denote the image of $g$ under our map as $h$}
    h\left(1\otimes \frac{1}{x^\alpha}\right)&= \frac{c_0^p}{x^{p-a} y^{\beta p- b}}\\
    h\left(1\otimes \frac{1}{x^\alpha+1}\right)&= \frac{c_0^p}{x^{2p-a}y^{\beta p-b}}+\frac{c_1^p}{x^{p-a}y^{(\beta+1)p-b}}\\
    &\vdots\\
    h\left(1\otimes \frac{1}{x^{\alpha+n}}\right)&=\sum\limits_{i=0}^{n} \frac{c_i^p}{x^{(n+1-i)p-a} y^{(\beta+i)p-b}}
    \end{align*}
Rewrite $h$ as $\sum h_i$ in a similar way (the indexing is changed as only terms which appear have powers  on $x$ congruent to $-a$ mod $p$ and $-b$ mod p for $y$)  to how we rewrote $f$ to define the map back to $F(\Hom(R_x,E))$.   

\[
    h_i\left( 1\otimes \frac{1}{x^{\alpha+n}}\right):=\frac{c_i^p}{x^{(n+1-i)p-a} y^{(\beta+i)p-b}}\ \ \ \  \forall n \geq i
    \]
    
    To compute $m_i,n_i, k_i, \ell_i$ we find that
 \[h_i\left( x^{p-a-1}y^{p-b-1}\otimes \frac{y^{\beta+i-1}}{x^{\alpha+i}}\right)= \frac{c_i^p}{xy}\]

So $m_i=p-a+1$, $n_i=p-b+1$, $k_i=(\beta+i-1)$, $\ell_i=(\alpha+i)$. Thus we compute $h$'s image going back to $F_*R\otimes {}^*\Hom(R_x,E)$
 \[h\mapsto x^{p-(p-a-1)-1}y^{p-(p-b-1)-1}\otimes \sum \hat{h}_i\]
 \begin{align*}
 \intertext{compute $\hat{h}_i$:}
 \hat{h}_i\left(\frac{1}{x^{\alpha+i+n}}\right)&=\frac{(c_i^p)^{\frac{1}{p}}}{x^{n+1}y^{\beta+i}}\ \ \ \ \ \forall n\geq 0\\
 \intertext{So we have}
  \hat{h}_i\left(\frac{1}{x^{\alpha+n}}\right)&=\frac{(c_i^p)^{\frac{1}{p}}}{x^{(n-i)+1}y^{\beta+i}}\ \ \ \ \ \forall n\geq i\\
 \intertext{Note for any $n\geq 0$ }
 \sum\limits_{i\geq 0} \hat{h}_i\left(\frac{1}{x^{\alpha+n}}\right)&=\sum\limits_{i=0}^{n}\hat{h}_i \left(\frac{1}{x^{\alpha+n}}\right)+\underbrace{\sum\limits_{i>n} \hat{h}_i\left(\frac{1}{x^{\alpha+n}}\right)}_{\text{all terms zero}}\\
 &=\sum\limits_{i=0}^{n} \frac{c_i}{x^{n-i+1}y^{\beta+i}}=g\left(\frac{1}{x^{\alpha+n}}\right)\\
\end{align*}
So we have $h\mapsto x^ay^b\otimes g$ as needed.
\\

Take $f\in {}^*\Hom(F_*R\otimes R_x, E)$ and write the ouptuts of $f$ in the way used when defining ${}^*\Hom (F_*R \otimes R_x, E)\longrightarrow F_*R\otimes{}^*\Hom(R_x,E) $.  Define $f_i$ and $\hat{f}_i$ in exactly the same way as well.  

\[
     f\mapsto \sum\limits_{0\leq r<p} \left( x^{p-m_i-1}y^{p-n_i-1}\otimes \sum\limits_{i\equiv r\  \text{mod}\ p} \hat{f}_i \right)\]
     Fix  $0\leq r<p$ and look at the image when mapped back to ${}^*\Hom(F_*R\otimes R_x, E)$:     
     
      \[\left( x^{p-m_i-1}y^{p-n_i-1}\otimes \sum\limits_{i\equiv r\  \text{mod}\ p} \hat{f}_i \right)\mapsto \phi \]
      \begin{align*}
      \phi\left(1\otimes\frac{1}{x^\ell}\right):&=  x^{p-m_i-1}y^{p-n_i-1} \sum\limits_{i\equiv r\  \text{mod}\ p} \hat{f}_i\left(\frac{1}{x^\ell}\right)^p\\
     &=  x^{p-m_i-1}y^{p-n_i-1} \sum\limits_{i\equiv r\  \text{mod}\ p} \frac{c_i}{x^{(\ell-\ell_i+1)p} y^{(k_i+1)p}}\\
      &=\sum\limits_{i\equiv r\  \text{mod}\ p}\frac{c_i}{x^{(\ell-\ell_i)p+m_i+1}y^{k_ip+1}}
      \end{align*}
      Note the $f_i$ are defined so that they only ever output one term when the input is a pure tensor of monomials, so
      \[f_i\left(x^{m_i}y^{n_i}\otimes \frac{y^{k_i}}{x^{\ell_i}}\right)= \frac{c_i}{xy}\implies f_i\left(1\otimes \frac{1}{x^\ell}\right)= \frac{c_i}{x^{(\ell-\ell_i)p+m_i+1}y^{k_ip+n_i+1}}\]
     So we have
      \[\left( x^{p-m_i-1}y^{p-n_i-1}\otimes \sum\limits_{i\equiv r\  \text{mod}\ p} \hat{f}_i \right)\mapsto \left( \sum\limits_{i\equiv r \text{ mod } p} f_i\right)\]
      So now just note that we can do this for each $r$ and then we get $\sum f_i=f$ as needed.
\end{proof}

\begin{remark}\label{rem1} $R_x$ has a graded $F$-module structure, given by the morphism $\theta_{R_x}(\frac{r}{x^n}):=rx^{n(p-1)}\otimes \frac{1}{x^n}$, with inverse ,$\theta_{R_x}^{-1}\left(r\otimes \frac{s}{x^n}\right):=\frac{rs^p}{x^{pn}}$. 
\end{remark}

\begin{prop}\label{prop3}$\sideset{^*}{_R}\Hom(R_x,E)$ has a graded $F$-module structure
\end{prop}
\begin{proof}
  First we get a degree preserving isomorphism \[\sideset{^*}{_R}\Hom(R_x,E)\xrightarrow{\sim}\sideset{^*}{_R}\Hom(F(R_x),E)\] by applying the above remark and Proposition \ref{prop1}.
  
We also have a degree preserving isomorphism \[\sideset{^*}{_R}\Hom(F(R_x),E)\xrightarrow{\sim} F(\sideset{^*}{_R}\Hom(R_x,E))\] constructed for Proposition \ref{prop2}.

Compose these two degree preserving isomorphisms to get a degree preserving isomorphism
\[\sideset{^*}{_R}\Hom(R_x,E)\xrightarrow{\sim}F(\sideset{^*}{_R}\Hom(R_x,E))\]
as needed.
\end{proof}

\begin{remark}\label{rem3}
By the same argument as in Theorem 4.5 of \cite{Ma14}, giving a graded $F$- module structure to $N$, where $N\cong R_x\oplus R$ as $R$-modules, is the same as giving a degree zero element of $R_x$.
\end{remark}

\begin{prop}\label{prop4}A graded $F$-Module structure on $R_x\oplus R$ induces a graded $F$-module structure on ${}^{*}\Hom(R_x,E)\oplus E$.
\end{prop}

\begin{proof}
Let $\theta_N$ be a graded $F$-module structure map for $N\cong R_x\oplus R$.  We need a degree preserving isomorphism ${}^*\Hom(R_x,E)\oplus E\rightarrow F({}^*\Hom(R_x,E)\oplus E)$.  We produce such a map by composing degree preserving isomorphisms

\[{}^*\Hom(R_x,E)\oplus E\rightarrow{}^*\Hom(R_x\oplus R,E)\rightarrow{}^*\Hom(F (R_x\oplus R),E)\rightarrow{}^*\Hom(F( R_x),E)\oplus {}^*\Hom(F(R),E)\]
\[{}^*\Hom(F( R_x),E)\oplus {}^*\Hom(F(R),E)\rightarrow F(\Hom( R_x,E))\oplus F(E)\rightarrow F({}^{*}\Hom(R_x,E)\oplus E)\]\\

We get  ${}^*\Hom(R_x,E)\oplus E\rightarrow{}^*\Hom(R_x\oplus R,E) $ via $(f,m)\mapsto g$ with $g(\frac{1}{x^n},0)=f(\frac{1}{x^n})$ and $g(0,1)=m$.

 Apply proposition \ref{prop1} to get a degree preserving isomorphism ${}^*\Hom(R_x\oplus R,E)\rightarrow {}^*\Hom(F (R_x\oplus R),E)$  induced by $\theta_N^{-1}: F (R_x\oplus R)\rightarrow R_x\oplus R$, i.e. $g\mapsto h$ where $h(r_1\otimes (\frac{r_2}{x^n},r_3))= g(\theta_M^{-1}(\frac{r_2}{x^n},r_3)))$.

 ${}^*\Hom(F (R_x\oplus R),E)\rightarrow {}^*\Hom(F( R_x),E)\oplus {}^*\Hom(F(R),E)$ is given by $h\mapsto (h_1,h_2)$ where $h_1(r_1\otimes \frac{r_2}{x^n})=h(r_1\otimes(\frac{r_2}{x^n},0))$ and $h_2(r_1\otimes r_2)=h(r_1\otimes(0,r_2))$. \\

To get ${}^*\Hom(F( R_x),E)\oplus {}^*\Hom(F(R),E\xrightarrow{\sim} F(\Hom( R_x,E))\oplus F(E)$ we  map the ${}^{*}\Hom(F(R),E)$ component to  $F(E)$ by $f\mapsto \sum_i c_i( x^{n_i}y^{m_i})^{p-1}\otimes \frac{1}{x^{n_i}y^{m_i}}$ where we have $f(1\otimes 1)= \sum_i\frac{c_i}{x^{n_i}y^{m_i}}$. Then we use the degree preserving isomorphism that we already constructed in Proposition \ref{prop2} to map ${}^{*}\Hom(F(R_x),E)$ to $F({}^{*}\Hom(R_x,E))$.

Finally $F(\Hom( R_x,E))\oplus F(E)\rightarrow F({}^{*}\Hom(R_x,E)\oplus E)$ is given by $(r_1\otimes f,r_2\otimes m)\mapsto r_1\otimes(f,0)+r_2\otimes(0,m)$.
\end{proof}

\section{Proof of Theorem \ref{main theorem}}
\label{proof of main theorem}

To prove Theorem \ref{main theorem} we construct a graded $F$-module, $L$, isomorphic to $\Hom(R_x,E)\oplus E$ as $R$-modules, such that the short exact sequence $E\rightarrow L \rightarrow  {}^*\Hom(R_x,E)$ does not split in $\mathcal{F}$ so $\Ext^{1}_{\mathcal{F}}({}^{*}\Hom(R_x,E),E)\neq 0$.

\begin{proof}[Proof of theorem 1.2]Define a graded $F$-module, $N$ where $N \cong R_x\oplus R$ as $R$-modules and, in the same manner as in Theorem 4.5 in \cite{Ma14}, we characterize the structure map $\theta_N$ by $\theta_{R_x}^{-1}\oplus \theta_R^{-1}\circ\theta_N(m,r)=(m+\frac{ry}{x},r)$.

Define $L$ such that $L\cong {}^*\Hom(R_x, E)\oplus E $  and $L$ has structure map induced by the structure map on $N$

Assume we have the following commutative diagram of graded $F$-modules

\begin{center}
\begin{tikzcd}
E \arrow[r] \arrow[d, "\text{id}_E"]   &{}^*\Hom(R_x, E)\oplus E \arrow[r] \arrow[d, "g"]& {}^*\Hom(R_x,E)\arrow[d, "\text{id}_{{}^*\Hom(R_x,E)}"]   \\
E\arrow[r] &  L \arrow[r] & {}^*\Hom(R_x,E)
\end{tikzcd}
\end{center}

Fix $g: {}^*\Hom(R_x, E)\oplus E\rightarrow L$  given by \[\left(\vp, \frac{r}{x^ky^m}\right)\mapsto \left(\vp, \frac{r}{x^ky^m}+ \vp\left(\frac{t}{x^\alpha}\right)\right)\]
for some homogeneous polynomial $t$ of degree $\alpha$ with a non-trivial term of $y$ degree $\alpha$ 



We show such a $g$ cannot exist by making an explicit computation for a particular element along different paths in the following commuting square

\begin{center}
\begin{tikzcd}
  {}^*\Hom(R_x, E)\oplus E \arrow[r,"g"] \arrow[d,"\theta_{(R_x)^\vee\oplus E}"] & L\arrow[d,"\theta_L"]   \\
    F( {}^*\Hom(R_x, E)\oplus E)\arrow[r,"\text{id}\otimes g"] & {}F(L)
\end{tikzcd}
\end{center}

Compute $\theta_L\circ g(\vp,\frac{1}{xy^{\alpha p+2}})$ 
where $\vp(\frac{1}{x^n}):= \frac{1}{x^{n+1}y^{\alpha p+2}}$

\[
g(\vp ,\frac{1}{xy^{\alpha p+2}})= (\vp, \frac{1}{xy^{\alpha p+2}}+\frac{t}{x^{\alpha+1}y^{\alpha p+2}})\]

Now we compute $\theta_L(\vp, \frac{1}{xy^{\alpha p+2}}+\frac{t}{x^{\alpha+1}y^{\alpha p+2}})$ according to the same steps as in Proposition \ref{prop4}.   

We first map $(\vp, \frac{1}{xy^{\alpha p+2}}+\frac{t}{x^{\alpha+1}y^{\alpha p+2}})$ to ${}^*\Hom(R_x\oplus R,E)$ to get $\vp'$ such that $\vp'(\frac{1}{x^n},0)=\vp(\frac{1}{x^n})$ and $\vp'(0,1)=\frac{1}{xy^{\alpha p+2}}+\frac{t}{x^{\alpha+1}y^{\alpha p+2}}$

 then apply the map induced by $\theta_N^{-1}$ from ${}^*\Hom(R_x\oplus R,E)$ to ${}^*\Hom(F (R_x\oplus R),E)$  we have $\psi$ where 
 \[\psi\left(1\otimes\left(\frac{1}{x^n},0\right)\right)=\vp'\left(\theta_N^{-1}\left(1\otimes\left(\frac{1}{x^n},0\right)\right)\right)=\frac{1}{x^{pn+1}y^{\alpha p+2}}\]
 \[\psi(1\otimes (0, 1))= \vp'\left(\frac{-y}{x},1\right)=\frac{-1}{x^2y^{\alpha p+1}}+\frac{1}{xy^{\alpha p+2}}+\frac{t}{x^{\alpha+1}y^{\alpha p+2}} \]
 
 If we then map $\psi$ to $\Hom(F(R_x),E)\oplus {}^{*}\Hom(F(R),E)$ then we have

\[
 \left(\psi_1, \psi_2\right)
 \]

where $\psi_1\in {}^{*}\Hom(F(R_x),E)$ is given by $\psi_1(1\otimes \frac{1}{x^n})=\psi\left(1\otimes \left(\frac{1}{x^n}, 0\right)\right)= \vp(\frac{1}{x^{pn}})= \frac{1}{x^{pn+1}y^{\alpha p+2}}$ and $\psi_2\in {}^*\Hom(F(R),E)$ is given by $\psi_2(1\otimes 1)=\psi(1\otimes (0,1))=\frac{-1}{x^2y^{\alpha p+1}}+\frac{1}{xy^{\alpha p+2}}+\frac{t}{x^{\alpha+1}y^{\alpha p+2}}$

When we map $(\psi_1,\psi_2)$ into $F({}^{*}\Hom(R_x,E))\oplus F(E)$ we get 

\[\left((xy^{\alpha p+2})^{p-1}\otimes \vp,-(x^2y^{\alpha p+1})^{p-1}\otimes \frac{1}{x^2y^{\alpha p+1}}+ (x^{\alpha+1}y^{\alpha p+2})^{p-1}t\otimes \frac{1}{x^{\alpha+1}y^{\alpha p+2}}+ (xy^{\alpha p+2})^{p-1}\otimes\frac{1}{xy^{\alpha p+2}}\right)\]

When we map the above line to $F ({}^*\Hom(R_x,E)\oplus E)$ we get
\begin{align}\label{eq1}
(xy^{\alpha p+2})^{p-1}\otimes(\vp,\frac{1}{xy^{\alpha p+2}})-(x^2y^{\alpha p+1})^{p-1}\otimes (0,\frac{1}{x^2y^{\alpha p+1}})+ (x^{\alpha+1}y^{\alpha p+2})^{p-1}t\otimes (0,\frac{1}{x^{\alpha+1}y^{\alpha p+2}})
\end{align}
Now we compare this result to $\text{id}\otimes g\circ \theta_{\Hom(R_x,E)\oplus E}(\vp,\frac{1}{xy^{\alpha p+2}})$
\begin{align*}
(\vp,\frac{1}{xy^{\alpha p+2}})&\mapsto (\psi,(xy^{\alpha p+2})^{p-1}\otimes \frac{1}{xy^{\alpha p+2}})\in {}^*\Hom(F(R_x),E)\oplus F(E)\\
\intertext{then  $(\psi,(xy^{\alpha p+2})^{p-1}\otimes \frac{1}{xy^{\alpha p+2}})$ maps to} 
 (xy^{\alpha p+2})^{p-1}\otimes (\vp,0)&+ (xy^{\alpha p+2})^{p-1}\otimes(0,\frac{1}{xy^{\alpha p+2}})\in F ({}^*\Hom(R_x,E)\oplus E)\\
 \end{align*}
 now applying $\text{id}\otimes g$
 \begin{align}\label{eq2}
 (xy^{\alpha p+2})^{p-1}\otimes (\vp,0)+ (xy^{\alpha p+2})^{p-1}\otimes (0,\frac{1}{xy^{\alpha p+2}})+(xy^{\alpha p+2})^{p-1}\otimes (0, \frac{t}{x^{\alpha+1}y^{\alpha p+2}})
\end{align}

Now apply $\theta_{{}^*\Hom(R_x,E)\oplus E}^{-1}$ to both \eqref{eq1} and \eqref{eq2} we have

\begin{align*}
    (\vp, \frac{1}{xy^{\alpha p+2}}-\frac{1}{x^2y^{\alpha p+1}}+\frac{t}{x^{\alpha +1}y^{\alpha p+2}})&=(\vp, \frac{1}{xy^{\alpha p+2}}+\frac{t^p}{x^{\alpha p+1}y^{\alpha p+2}})
    \intertext{so we have}
    \frac{-1}{x^2y^{\alpha p+1}}+\frac{t}{x^{\alpha+1}y^{\alpha p +2}}&= \frac{t^p}{x^{\alpha p+1}y^{\alpha p+2}}\\
    \intertext{If $\alpha=0$ then $t$ is a scalar and we would have}
    \frac{t^p x-tx+y}{x^2y^2}=0\\
    \intertext{since $x^2,y^2$ is regular in $R$, we would have}
(t^p-t)x+y&\in (x^{2},y^{2})\\
\intertext{which is impossibble because of degree considerations}
\end{align*}

If $\alpha>0$ then 
 \[\frac{t^p-t(x^{\alpha(p-1)})+yx^{\alpha p-1}}{x^{\alpha p+1}y^{\alpha p +2}}=0 \]
Note again that $x^{\alpha p+1}, y^{\alpha p+2}$ is regular, thus
 \[
    \implies t^p-tx^{\alpha(p-1)}+yx^{\alpha p-1} \in (x^{\alpha p+1},y^{\alpha p+2})\]
  Recall that $t$ is a homogeneous polynomial of degree $\alpha$, with non-trivial term of $y$ degree $\alpha$.  Thus we may write $t$ explicitly as $t=b_0y^\alpha +b_1y^{\alpha-1}x+\cdots+b_\alpha x^{\alpha}$ for some scalars $b_0,\ldots b_\alpha$. Then we would have

     \[
    \implies b_0^py^{\alpha p}+b_1^px^p y^{(\alpha -1)p}+\ldots +b_{\alpha}^p x^{\alpha p} - (b_0 y^\alpha x^{\alpha (p-1)}+\ldots b_\alpha x^{\alpha p})+(+yx^{\alpha p-1}) \in (x^{\alpha p+1},y^{\alpha p+2})\]
    Now multiply the whole lefthand side by $x^{\alpha p}$ to get
    \[\implies b_0^p x^{\alpha p}y^{\alpha p}+ \underbrace{x^{\alpha p}\left(b_1^px^p y^{(\alpha -1)p}+\ldots +b_{\alpha}^p x^{\alpha p} - (b_0 y^\alpha x^{\alpha (p-1)}+\ldots b_\alpha x^{\alpha p})+(+yx^{\alpha p-1})\right)}_{\in (x^{\alpha p+1})\subset (x^{\alpha p+1},y^{\alpha p+2}) } \in (x^{\alpha p+1},y^{\alpha p+2})\]
   $ \implies b_0x^{\alpha p}y^{\alpha p}\in (x^{\alpha p+1},y^{\alpha p+2})$
   
   Which is again impossible due to degree considerations.

Thus we have no such commuting diagram and by Yoneda's characterization of $\Ext^1$ we have that $\Ext_{\mathcal{F}}^{1}\left((R_x)^\vee,E\right)\neq 0$.
\end{proof}


\end{document}